\documentclass[11pt]{article}
\usepackage{enumerate}
\usepackage{amssymb,a4wide,latexsym,makeidx,epsfig,fleqn}
\usepackage{amsthm}
\usepackage{amsmath}
\usepackage{enumerate}
\usepackage{mathrsfs}

\newtheorem{theorem}{Theorem}[section]

\newtheorem{definition}[theorem]{Definition}
\newtheorem{lemma}[theorem]{Lemma}
\newtheorem{fact}[theorem]{Fact}

\newtheorem{claim}[theorem]{\hspace{6mm}Claim}
\newtheorem{proposition}[theorem]{Proposition}
\newtheorem{corollary}[theorem]{Corollary}

\newtheorem{question}[theorem]{Question}
\newtheorem{conjecture}[theorem]{Conjecture}

\begin{document}
\textwidth 150mm \textheight 225mm
\title{Forbidden rainbow subgraphs that force large monochromatic or multicolored $k$-connected subgraphs
\thanks{Supported by the National Natural Science Foundation of China (No. 11871398), the Natural Science Basic Research Plan in Shaanxi Province of China (Program No. 2018JM1032), and the Fundamental Research Funds for the Central Universities (No. 3102019ghjd003).}}
\author{{Xihe Li$^{a,b}$ and Ligong Wang$^{a,b,}$\thanks{Corresponding author.}}\\
{\small $^a$Department of Applied Mathematics, School of Science,}\\ {\small Northwestern Polytechnical University, Xi'an, Shaanxi 710072, P. R. China. }\\ {\small $^b$Xi'an-Budapest Joint Research Center for Combinatorics,}\\ {\small Northwestern Polytechnical University, Xi'an, Shaanxi 710129, P. R. China. }\\ {\small E-mail: lxhdhr@163.com; lgwangmath@163.com}}
\date{}
\maketitle
\begin{center}
\begin{minipage}{120mm}
\vskip 0.3cm
\begin{center}
{\small {\bf Abstract}}
\end{center}
{\small We consider a forbidden subgraph condition that implies the existence of a large highly connected monochromatic or multicolored subgraph. Let $n, k, m$ be positive integers with $n\gg m\gg k$, and let $\mathcal{A}$ be the set of graphs $G$ of order at least 3 such that there is a $k$-connected monochromatic subgraph of order at least $n-f(G,k,m)$ in any rainbow $G$-free coloring of $K_n$ using all the $m$ colors. In this paper, we prove that the set $\mathcal{A}$ consists of precisely $P_6$, $P_3\cup P_4$, $K_2\cup P_5$, $K_2\cup 2P_3$, $2K_2\cup K_3$, $2K_2\cup P^{+}_4$, $3K_2\cup K_{1,3}$ and their subgraphs of order at least 3. Moreover, we show that for any graph $H\in \mathcal{A}$, if $n$ sufficiently larger than $m$ and $k$, then any rainbow $(P_3\cup H)$-free coloring of $K_n$ using all the $m$ colors contains a $k$-connected monochromatic subgraph of order at least $cn$, where $c=c(H)$ is a constant, not depending on $n$, $m$ or $k$.

\ \ \ Furthermore, we consider a parallel problem in complete bipartite graphs. Let $s, t, k, m$ be positive integers with ${\rm min}\left\{s, t\right\}\gg m\gg k$ and $m\geq |E(H)|$, and let $\mathcal{B}$ be the set of bipartite graphs $H$ of order at least 3 such that there is a $k$-connected monochromatic subgraph of order at least $s+t-f(H,k,m)$ in any rainbow $H$-free coloring of $K_{s,t}$ using all the $m$ colors, where $f(H,k,m)$ is not depending on $s$ or $t$. We prove that the set $\mathcal{B}$ consists of precisely $2P_3$, $2K_2\cup K_{1,3}$ and their subgraphs of order at least 3.

\ \ \ Finally, we consider the large $k$-connected multicolored subgraph instead of monochromatic subgraph. We show that for $1\leq k \leq 3$ and $n$ sufficiently large, every Gallai-3-coloring of $K_n$ contains a $k$-connected subgraph of order at least $n-\left\lfloor\frac{k-1}{2}\right\rfloor$ using at most two colors. We also show that the above statement is false for $k=4t$, where $t$ is an positive integer.

 }
\end{minipage}
\end{center}

\section{Introduction }
\label{sec:ch-introduction}

Throughout this article, we only consider edge-colorings of finite simple graphs. Given a graph $G=(V(G),$ $E(G))$ and an integer $m \geq 1$, let $c$ : $E(G)\rightarrow [m]$ be a {\it $m$-coloring} on the edges of $G$, where $[m]= \{1, 2, \ldots, m\}$ is the set of colors. For convenient, we also use names of colors like ``red", ``blue" or ``green". A coloring of a graph is called {\it rainbow} if all the edges are colored differently, and {\it monochromatic} if all the edges have the same color. For two graphs $G$ and $H$, let $G\vee H$ denote the join of $G$ and $H$, $G\cup H$ denote the union of $G$ and $H$, and $nG$ denote the union of $n$ disjoint copies of $G$.

It is easy to check that any 2-coloring of $K_n$ contains a 1-connected monochromatic subgraph, a spanning tree. In general, it is not possible to find a spanning monochromatic subgraph with higher connectivity for every 2-colored $K_n$. For instance, consider a 2-coloring of $K_n$ using red and blue, in which the monochromatic subgraph induced by red edges is a spanning star (only 1-connected). Then the spanning monochromatic subgraph consisting of all the blue edges is $K_1\cup K_{n-1}$ that is not connected. Thus if we would like to find a monochromatic subgraph with higher connectivity, then we must consider the subgraph with order smaller than $K_n$. In 2008, Bollob\'{a}s and Gy\'{a}rf\'{a}s \cite{BoGy} posed the following conjecture, in which they considered the monochromatic subgraph of order at least $n-2(k-1)$.

\noindent\begin{conjecture}\label{conj-1} {\normalfont (\cite{BoGy})}
If $n>4(k-1)$, then every 2-coloring of $K_n$ contains a monochromatic $k$-connected subgraph of order at least $n-2(k-1)$.
\end{conjecture}

In \cite{FuMa2}, Fujita and Magnant proved that Conjecture \ref{conj-1} holds for $n>6.5(k-1)$.

\noindent\begin{theorem}\label{th:2coloring} {\normalfont (\cite{FuMa2})}
If $n>6.5(k-1)$, then every 2-coloring of $K_n$ contains a monochromatic $k$-connected subgraph of order at least $n-2(k-1)$.
\end{theorem}

Moreover, for more colors, we shall consider the {\it almost spanning $k$-connected monochromatic subgraph} defined as follows.

\noindent\begin{definition}\label{de:ASMS} {\normalfont (\cite{FuMa3})}
Given positive integers $k$ and $m$, for sufficiently large $n$ compared with $k$ and $m$, an almost spanning $k$-connected monochromatic subgraph, denoted by $ASMS(k)$, is a $k$-connected monochromatic subgraph of order at least $n-f(k,m)$ in a colored $K_n$ using all the $m$ colors, where $f(k,m)$ is not depending on $n$.
\end{definition}

When we consider a coloring of $K_n$ using $m$ colors, Liu, Morris and Prince \cite{LiMP2} showed that the best order for a highly connected monochromatic subgraph would be $\frac{n}{m-1}$. Thus, to find an $ASMS(k)$, we must consider more additional restrictions on the coloring of $K_n$, such as forbidden proper subgraphs \cite{KaMN} and forbidden rainbow subgraphs \cite{FuMa3}. Let $P^{+}_4$ be the tree with degree sequence 1, 1, 1, 2, 3. In \cite{FuMa3}, Fujita and Magnant proved the following result.

\noindent\begin{theorem}\label{th-1} {\normalfont (\cite{FuMa3})}
Let $n, k, m$ be positive integers with $n\gg m\gg k$ and let $\mathcal{G}$ be the set of connected graphs $G$ of order at least 3 such that there is an $ASMS(k)$ in any rainbow $G$-free coloring of $K_n$ using all the $m$ colors. Then $\mathcal{G}$ consists of precisely $K_3$, $P_6$, $P^{+}_4$ and their subgraphs of order at least 3.
\end{theorem}

In this paper, we first consider the following question.

\noindent\begin{question}\label{ques-1}
For what disconnected graphs $H$ of order at least 3 without isolated vertex does the following statement holds? Let $n, k, m$ be positive integers with $n\gg m\gg k$ and $m\geq |E(H)|$. There is an $ASMS(k)$ in every rainbow $H$-free coloring of $K_n$ using all the $m$ colors.
\end{question}

Note that since we only consider edge-colorings, we will not consider the graphs with isolated vertices. Let $\mathcal{H}$ be the set of disconnected graphs $H$ such that $H$ satisfies Question \ref{ques-1}. For a graph $G$, the {\it component number} $c_0(G)$ is the number of components of $G$. For integers $i\geq 2$, let $\mathcal{H}^{(i)}=\{H\ |\ c_0(H)=i, H\in \mathcal{H}\}$ and define $\mathcal{H}^{(1)}=\mathcal{G}$. Clearly, we have $\mathcal{H}=\bigcup_{i\geq 2}\mathcal{H}^{(i)}$. We have the following main result.

\noindent\begin{theorem}\label{th:discon}
The set $\mathcal{H}$ consists of precisely $P_3\cup P_4$, $K_2\cup P_5$, $K_2\cup 2P_3$, $2K_2\cup K_3$, $2K_2\cup P^{+}_4$, $3K_2\cup K_{1,3}$ and their disconnected subgraphs of order at least 3.
\end{theorem}

Let $\mathcal{A}=\mathcal{G}\cup \mathcal{H}$. Then clearly we have that $\mathcal{A}$ consists of precisely $P_6$, $P_3\cup P_4$, $K_2\cup P_5$, $K_2\cup 2P_3$, $2K_2\cup K_3$, $2K_2\cup P^{+}_4$, $3K_2\cup K_{1,3}$ and their subgraphs of order at least 3.

For other rainbow subgraph $G\notin \mathcal{G}\cup \mathcal{H}$, the largest monochromatic $k$-connected subgraph will not be almost spanning, but how large should it be? We prove the following result.

\noindent\begin{theorem}\label{th:main2}
For any graph $H\in \mathcal{G}\cup \mathcal{H}$, if $n$ sufficiently larger than $m$ and $k$, then any rainbow $(P_3\cup H)$-free coloring of $K_n$ using all the $m$ colors contains a $k$-connected monochromatic subgraph of order at least $cn$, where $c=c(H)$ is a constant, not depending on $n$, $m$ or $k$.
\end{theorem}

Another natural question is to consider a parallel problem in complete bipartite graphs instead of complete graphs.

\noindent\begin{question}\label{ques-2}
For what bipartite graphs $H$ of order at least 3 without isolated vertex does the following statement holds? Let $s, t, k, m$ be positive integers with ${\rm min}\left\{s, t\right\}\gg m\gg k$ and $m\geq |E(H)|$. There is a $k$-connected monochromatic subgraph of order at least $s+t-f(H,k,m)$ in every rainbow $H$-free coloring of $K_{s,t}$ using all the $m$ colors, where $f(H,k,m)$ is not depending on $s$ and $t$.
\end{question}

Let $\mathcal{B}$ be the set of bipartite graphs $H$ satisfying Question \ref{ques-2}. We have the following result.

\noindent\begin{theorem}\label{th:bipartite}
The set $\mathcal{B}$ consists of precisely $2P_3$, $2K_2\cup K_{1,3}$ and their subgraphs of order at least 3.
\end{theorem}

Moreover, instead of looking for monochromatic subgraphs, we try to find large 2-colored subgraphs in Gallai-3-colorings of complete graphs. Here a {\it Gallai-$m$-coloring} means a rainbow triangle-free coloring with $m$ colors appearing. In \cite{FoGP}, Fox, Grinshpun and Pach proved that every Gallai-3-coloring of $K_n$ contains a complete graph of order $\Omega(n^{1/3}\rm{log}^2n)$ that uses at most two colors. In \cite{Wag}, Wagner proved that every Gallai-3-coloring of $K_n$ contains a 2-colored subgraph with chromatic number at least $n^{2/3}$. In this paper, we consider the large $k$-connected subgraphs using at most two colors in Gallai-3-colorings. By the following example we know that the largest $k$-connected subgraphs using at most two colors has order at most $n-\left\lfloor\frac{k-1}{2}\right\rfloor$. Let $G$ be a $3$-coloring of $K_n$ with $V(G)=V_1\cup V_2\cup V_3$, such that $|V_1|=n-k+1$, $|V_2|=\left\lceil\frac{k-1}{2}\right\rceil$ and $|V_3|=\left\lfloor\frac{k-1}{2}\right\rfloor$. We color the edges of $G$ such that $c(V_1, V_2)=c(V_2, V_3)=1$, $c(V_1, V_3)=2$ and $c(V_1)=c(V_2)=c(V_3)=3$. It is easy to see that $G$ is a Gallai-3-coloring of $K_n$, and the largest $k$-connected subgraph using at most two colors has order $n-\left\lfloor\frac{k-1}{2}\right\rfloor$. In light of this example, we may consider the following question.

\noindent\begin{question}\label{ques-3}
For $k \geq 1$ and $n$ sufficient large, does every Gallai-3-coloring of $K_n$ contains a $k$-connected subgraph of order at least $n-\left\lfloor\frac{k-1}{2}\right\rfloor$ using at most two colors?
\end{question}

We will prove that Question \ref{ques-3} is true for $1\leq k \leq 3$, and we will give a counterexample to show that this question is false when $k= 4t$, where $t$ is an positive integer.

Finally, we define some notations and terminologies. For $U\subseteq V(G)$, let $G[U]$ denote the subgraph of $G$ induced by $U$, and let $G-U$ denote the graph obtained from $G$ by deleting all the vertices of $U$ along with the edges incident with some vertex in $U$. For two graphs $G$ and $H$, let $G\subseteq H$ denote that $G$ is a subgraph of $H$. For an edge $uv$, let $c(uv)$ denote the color used on the edge $uv$. For $U, V \subseteq V(G)$, let $C(U, V)$ (resp., $C(U)$) be the set of colors used on the edges between $U$ and $V$  (resp., within $U$). Given a graph $G$, we use $\Delta(G)$ to denote the maximum degree of $G$.

The remainder of this paper is organized as follows. In Section \ref{sec:ch2}, we give the complete characterisation of $\mathcal{H}$. In Section \ref{sec:ch3}, we give the proof of Theorem \ref{th:main2}. In Section \ref{sec:ch4*}, we give a characterisation of $\mathcal{B}$. In Section \ref{sec:ch4}, we will consider Question \ref{ques-3}. In Section \ref{sec:ch5}, we conclude this paper by presenting some related future works.

\section{A classification of the set $\mathcal{H}$}
\label{sec:ch2}

We begin with the following lemmas and corollary which will be used in the proof of our main result.

\noindent\begin{lemma}\label{le:conK3} {\normalfont (\cite{FuMa3})}
Suppose $m\geq 3$ and $k\geq 2$. If $n\geq (m+11)(k-1)+7k {\rm ln} k$, then every Gallai-$m$-coloring of $K_n$ contains a $k$-connected monochromatic subgraph of order at least $n-m(k-1)$.
\end{lemma}

\noindent\begin{lemma}\label{le:conK13P4+} {\normalfont (\cite{FuMa3})}
For $G\in \{K_{1,3}, P^{+}_4\}$, if $m\geq 4$, $k\geq 1$ and $n\geq 7.5(k-1)$, then every rainbow $G$-free coloring of $K_n$ using all the $m$ colors contains a $k$-connected monochromatic subgraph of order at least $n-3k+2$.
\end{lemma}

\noindent\begin{lemma}\label{le:conP6} {\normalfont (\cite{FuMa3})}
Suppose $k\geq 1$ and $m\geq {\rm max}\left\{\frac{k}{2}+8, 15\right\}$. If $n\geq (m+11)(k-1)+7k {\rm ln} k+2m+3$, then every rainbow $P_6$-free coloring of $K_n$ using all the $m$ colors contains a $k$-connected monochromatic subgraph of order at least $n-7k+2$.
\end{lemma}

\noindent\begin{corollary}\label{co:K2k3K2P4+K2P5P3P4}
Suppose $m\geq {\rm max}\left\{\frac{k}{2}+8, 44\right\}$ and $k\geq 2$. For every $H\in \{K_2\cup K_3, K_2\cup P^+_4, K_2\cup P_5, P_3\cup P_4\}$, if $n\geq (m+11)(k-1)+7k {\rm ln} k+2m+3$, then every rainbow $H$-free coloring of $K_n$ using all the $m$ colors contains an $ASMS(k)$ {\rm(}that is, $H \in \mathcal{H}^{(2)}${\rm)}.
\end{corollary}

{\bf Proof.} For any $H\in \{K_2\cup K_3, K_2\cup P^+_4, K_2\cup P_5, P_3\cup P_4\}$, let $G$ be a coloring of $K_n$ using all the $m$ colors without rainbow $H$, and suppose that there exists no $ASMS(k)$.

Firstly, we consider the case when $H=K_2\cup K_3$. By Lemma \ref{le:conK3}, there is a rainbow $K_3$ in $G$, say with vertex set $V=\{v_1, v_2, v_3\}$ such that $c(v_1v_2)=1$, $c(v_2v_3)=2$ and $c(v_3v_1)=3$. Then $C(V(G)\setminus V)\subseteq \{1, 2, 3\}$ to avoid a rainbow $K_2\cup K_3$. If $|C(V(G)\setminus V)|\leq 2$, then there is an $ASMS(k)$ by Theorem \ref{th:2coloring}. Thus $C(V(G)\setminus V)=\{1, 2, 3\}$. Since $|V(G)\setminus V|\geq (m+11)(k-1)+7k {\rm ln} k+2m$, there is also a rainbow $K_3$ in $G[V(G)\setminus V]$ by Lemma \ref{le:conK3}, say with vertex set $V'=\{v_4, v_5, v_6\}$. Since there are $m\geq {\rm max}\left\{\frac{k}{2}+8, 44\right\}>12$ colors used in $G$, there is an edge $e$ such that $c(e)\geq 4$ and $e$ is incident with some vertex in $V(G)\setminus (V\cup V')$. Then there is a rainbow $K_2\cup K_3$, a contradiction.

Secondly, we consider the case when $H=K_2\cup P^+_4$. By Lemma \ref{le:conK13P4+}, there is a rainbow $P^+_4$ in $G$, say with vertex set $U=\{u_1, u_2, \ldots, u_5\}$ such that $c(u_1u_2)=1$, $c(u_2u_3)=2$, $c(u_3u_4)=3$ and $c(u_3u_5)=4$. Then $C(V(G)\setminus U)\subseteq \{1, 2, 3, 4\}$ to avoid a rainbow $K_2\cup P^+_4$. If $|C(V(G)\setminus U)|\leq 2$ , then there is an $ASMS(k)$ by Theorem \ref{th:2coloring}, a contradiction. If $|C(V(G)\setminus U)|=4$, then by Lemma \ref{le:conK13P4+}, there is also a rainbow $P^+_4$ in $G[V(G)\setminus U]$, say with vertex set $U'=\{u_6, u_7, \ldots, u_{10}\}$. Since there are $m\geq 44$ colors used in $G$, there is an edge $e$ such that $c(e)\geq 5$ and $e$ is incident with some vertex in $V(G)\setminus (U\cup U')$. Then there is a rainbow $K_2\cup P^+_4$, a contradiction. Thus $|C(V(G)\setminus U)|=3$, say $C(V(G)\setminus U)=\{c_1, c_2, c_3\}$  and $c(e_i)=c_i$ for $i=1, 2, 3$, where $e_1$, $e_2$ and $e_3$ are three edges in $G[V(G)\setminus U]$. Let $U''$ be the set of vertices incident with at least one of $e_1$, $e_2$ and $e_3$, so $|U''|\leq 6$. Since there are $m\geq 44$ colors used in $G$, there are at least $m-\left(|U| |U''|+\binom{|U|}{2}\right)\geq 4$ distinct colors, say $c'_1, c'_2, c'_3$ and $c'_4$, on the edges between $U$ and $V(G)\setminus (U\cup U''))$ with $\{c'_1, c'_2, c'_3, c'_4\}\subseteq \{5, 6, \ldots, m\}$. It is easy to check that there is a rainbow $K_2\cup P^+_4$, a contradiction.

Finally, we consider the case when $H=K_2\cup P_5$ or $H=P_3\cup P_4$. By Lemma \ref{le:conP6}, there is a rainbow $P_6$ in $G$, say with vertex set $W=\{w_1, w_2, \ldots, w_6\}$ such that $c(w_iw_{i+1})=i$ for every $1\leq i\leq 5$. When $H=K_2\cup P_5$, we have $C(V(G)\setminus W)\subseteq \{2, 3, 4\}$ for avoiding a rainbow $K_2\cup P_5$, and moreover we have $C(\{w_1, w_3, w_4, w_6\}, V(G)\setminus W)\subseteq \{1, 2, \ldots, 5\}$ for the same reason. Thus we may assume that $c(w_2w')=6$ for some $w'\in V(G)\setminus W$. Then $C(w', V(G)\setminus (W\cup\{w'\}))\subseteq \{2, 3\}$. If $C(w', V(G)\setminus (W\cup\{w'\}))= \{2, 3\}$, say $c(w'w'')=2$ and $c(w'w''')=3$, where $w''$ and $w'''$ are two distinct vertices in $V(G)\setminus (W\cup\{w'\})$, then $c(w''w''')\notin \{1, 2, \ldots, m\}$, a contradiction. If $c(w', V(G)\setminus (W\cup\{w'\}))=2$ (resp., $c(w', V(G) \setminus (W\cup\{w'\}))=3$), then $C(V(G)\setminus (W\cup\{w'\})=2$ (resp., $C(V(G)\setminus (W\cup\{w'\})=3$), resulting in an $ASMS(k)$, a contradiction. When $H=P_3\cup P_4$, we have $C(W, V(G)\setminus W)\subseteq \{1, 2, \ldots, 5\}$ for avoiding a rainbow $P_3\cup P_4$. Thus we may assume that $c(xy)=6$, where $x$ and $y$ are two distinct vertices in $V(G)\setminus W$. If $c(w_2x) =1$, then $w_2xy$ and $w_3w_4w_5w_6$ form a rainbow $P_3 \cup P_4$, a contradiction. If $c(w_2x)=5$ (resp., $c(w_2x)\in \{2, 3\}$), then $w_3w_4w_5$ (resp., $w_4w_5w_6$) and $w_1w_2xy$ form a rainbow $P_3 \cup P_4$, a contradiction. Thus $c(w_2x)=4$, and by symmetry we have $c(w_2y)=4$, $c(w_5x)=c(w_5y)=2$. For any vertex $z\in V(G)\setminus (W\cup \{x, y\})$. If $c(xz)\notin \{3, 6\}$, then it is easy to find a rainbow $P_3 \cup P_4$. If $c(xz)\in \{3, 6\}$, then $w_1w_2xz$ and $yw_5w_6$ form a rainbow $P_3\cup P_4$, a contradiction.
\hfill$\blacksquare$
\vspace{0.1cm}

Note that if $G$ and $H$ be two graphs with $G\subseteq H$, then any rainbow $G$-free coloring of $K_n$ certainly contains no rainbow $H$. Thus we have the following simple result.

\noindent\begin{fact}\label{fa:2}
Let $G$ and $H$ be two graphs with $G\subseteq H$ and $|V(G)|\geq 3$. If $G \notin \mathcal{G}\cup \mathcal{H}$, then $H\notin \mathcal{H}$.
\end{fact}

Now we give a complete characterisation of $\mathcal{H}^{(2)}$.

\noindent\begin{theorem}\label{th:H2}
The set $\mathcal{H}^{(2)}$ consists of precisely $P_3\cup P_4$, $K_2\cup P_5$, $K_2\cup K_3$, $K_2\cup P^{+}_4$ and their subgraphs of order at least 3 with component number 2.
\end{theorem}

\noindent {\bf Proof.} In the following, we will prove that $\mathcal{H}^{(2)}$ is a subset of $P_3\cup P_4$, $K_2\cup K_3$, $K_2\cup P_5$, $K_2\cup P^{+}_4$ and their subgraphs of order at least 3 with component number 2, which together with Corollary \ref{co:K2k3K2P4+K2P5P3P4} and Fact \ref{fa:2} completes the proof of Theorem \ref{th:H2}.

For the proof, we will consider two colorings $R_1$ and $R_2$ (see Fig. 2.1) constructed as follows. Let $R$ be an $m$-coloring of $K_n$ with $V(R)=V_1\cup V_2\cup V_3$, such that each of $V_1$, $V_2$ and $V_3$ contains about $\frac{n}{3}$ vertices. We color the edges such that $C(V_1, V_2)=\{1\}$, $C(V_2)=C(V_2, V_3)=\{2\}$, $C(V_3)=C(V_3, V_1)=\{3\}$, and all the edges within $V_1$ are colored with color 1 except for a rainbow matching (resp., rainbow star) using all the remaining colors, and let $R_1=R$ (resp., $R_2=R$). Note that both the largest monochromatic $k$-connected subgraphs of $R_1$ and $R_2$ have order about $\frac{2n}{3}$, which are not ASMS($k$)s. Thus every graph in $\mathcal{H}^{(2)}$ should be a rainbow subgraph of both $R_1$ and $R_2$..

\begin{figure}[htb]
\begin{center}
\includegraphics {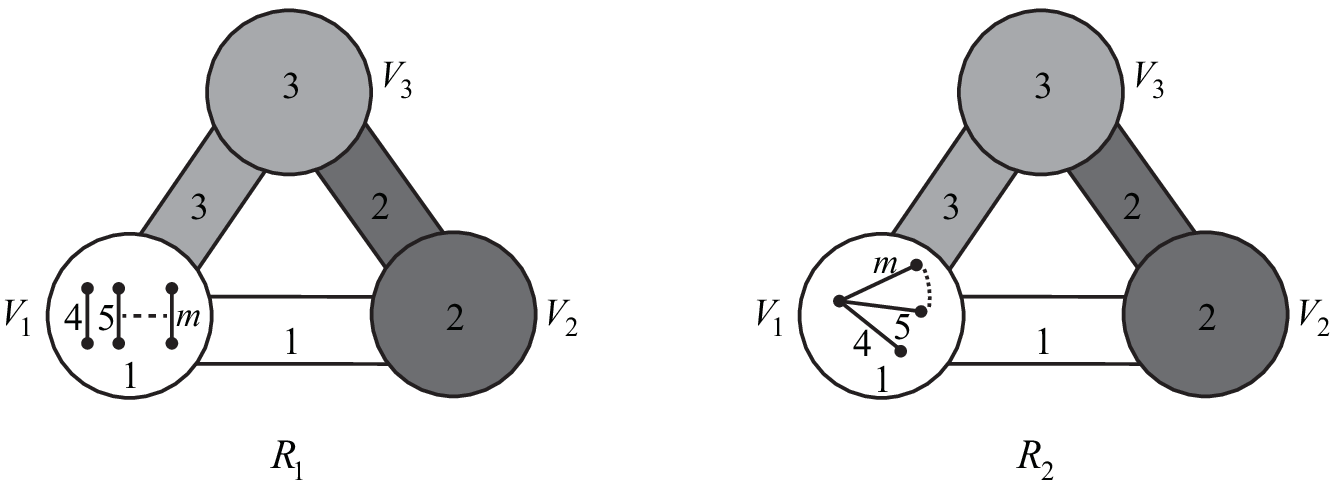}
\centerline{Fig. 2.1: The colorings $R_1$ and $R_2$.}
\end{center}
\end{figure}

Recall that we have $\mathcal{G}=\{K_3, P_3, P_4, P_5, P_6, K_{1,3}, P^{+}_4\}$ (see Theorem \ref{th-1}). By Fact \ref{fa:2}, each graph in $\mathcal{H}$ is a union of graphs in $\mathcal{G}\cup\{K_2\}$. In the rest of the proof, we consider six types of graphs that are possible in $\mathcal{H}^{(2)}$ (note that some graphs might be contained in two types).

Type 1. $K_2\cup G$, where $G\in \mathcal{G}$;

Type 2. $P_6\cup G$, where $G\in \mathcal{G}$;

Type 3. $K_3\cup G$, where $G\in \mathcal{G}$;

Type 4. $P_5\cup G$, where $G\in \mathcal{G}$;

Type 5. $K_{1,3}\cup G$ and $P^+_4\cup G$, where $G\in \mathcal{G}$;

Type 6. $P_3\cup G$ and $P_4\cup G$, where $G\in \mathcal{G}$.

Since every rainbow $K_2\cup P_6$ contains six distinct colors (at least three colors in $\{4, 5, \ldots,$ $m\}$), there is no rainbow $K_2\cup P_6$ in $R_2$. Thus $K_2\cup P_6 \notin \mathcal{H}^{(2)}$ and there is no graph of Type 2 that is possible in $\mathcal{H}^{(2)}$. Moreover, it is easy to find a rainbow $K_2\cup K_3$, $K_2 \cup P_5$ and $K_2 \cup P^+_4$ in both $R_1$ and $R_2$. Thus $K_2\cup K_3$, $K_2 \cup P_5$, $K_2 \cup P^+_4$ along with their subgraphs of order at least 3 and component number 2 are possible in $\mathcal{H}^{(2)}$.

Since every rainbow $K_3$ is colored by colors 1, 2 and 3, every rainbow $K_{1,3}$ contains colors 1 and 3, every rainbow $P_5$ contains colors 1 and 3, and every rainbow $P_3$ contains an edge with color 1 or color 3 in $R_1$, we have $K_3 \cup P_3, K_{1,3} \cup P_3, P^+_4 \cup P_3, P_5 \cup P_3 \notin \mathcal{H}^{(2)}$. Thus there is no graph of Types 3, 4 and 5 that is possible in $\mathcal{H}^{(2)}$. For the same reason with $K_2\cup P_6$, we have $2P_4 \notin \mathcal{H}^{(2)}$. Thus only $2P_3$ and $P_3\cup P_4$ of Type 6 are possible in $\mathcal{H}^{(2)}$.
\hfill$\blacksquare$
\vspace{0.1cm}

In the following, we will consider the set $\mathcal{H}$.

\noindent\begin{lemma}\label{le:H}
The set $\mathcal{H}$ is a subset of $P_3\cup P_4$, $K_2\cup P_5$, $K_2\cup 2P_3$, $2K_2\cup K_3$, $2K_2\cup P^{+}_4$, $3K_2\cup K_{1,3}$ and their disconnected subgraphs of order at least 3.
\end{lemma}

\noindent {\bf Proof.} By Fact \ref{fa:2}, each graph in $\mathcal{H}$ is a union of graphs in $\mathcal{G}\cup\{K_2\}$. Thus for every $H\in \mathcal{H}$, we have $\Delta(H)\leq 3$. Using the coloring $R_2$ constructed in the proof of Theorem \ref{th:H2}, we can deduce that
\begin{center}
$H$ has at most four components,
\end{center}
and
\begin{center}
if $\Delta(H)=1$, then $e(H)\leq 4$; if $\Delta(H)=2$, then $e(H)\leq 5$; if $\Delta(H)=3$, then $e(H)\leq 6$.
\end{center}
Moreover, by Fact \ref{fa:2} we have that for every $H\in \mathcal{H}$,
\begin{center}
the subgraph induced by any two components of $H$ is a graph in $\mathcal{H}^{(2)}$.
\end{center}
Now for every graph $H'\in \mathcal{H}^{(2)}=\{P_3\cup P_4, 2P_3, K_2\cup P_4, K_2\cup P_3, 2K_2, K_2\cup P_5, K_2\cup K_3, K_2\cup P^+_4, K_2\cup K_{1,3}\}$, we can construct graphs that are possible in $\mathcal{H}$ by adding components to $H'$ subject to the above properties. The proof is complete.
\hfill$\blacksquare$

\noindent\begin{lemma}\label{le:H+}
Suppose $m\geq {\rm max}\left\{\frac{k}{2}+8, 77\right\}$ and $k\geq 2$. For every $H\in \{K_2\cup 2P_3$, $2K_2\cup K_3$, $2K_2\cup P^{+}_4$, $3K_2\cup K_{1,3}\}$, if $n\geq (m+11)(k-1)+7k {\rm ln} k+2m+3$, then every rainbow $H$-free coloring of $K_n$ using all the $m$ colors contains an $ASMS(k)$ {\rm(}that is, $H \in \mathcal{H}${\rm)}.
\end{lemma}

{\bf Proof.} For any $H\in \{K_2\cup 2P_3$, $2K_2\cup K_3$, $2K_2\cup P^{+}_4$, $3K_2\cup K_{1,3}\}$, let $G$ be a rainbow $H$-free coloring of $K_n$ using all the $m$ colors, and suppose that there exists no $ASMS(k)$.

Firstly, we consider the case when $H=K_2\cup 2P_3$. By Lemma \ref{le:conP6}, there is a rainbow $P_6$ in $G$, say with vertex set $U=\{u_1, u_2, \ldots, u_6\}$ such that $c(u_iu_{i+1})=i$ for every $1\leq i\leq 5$. For avoiding a rainbow $K_2\cup 2P_3$, we have $C(V(G)\setminus U)\subseteq \{1, 2, 4, 5\}$ and $|C(u_i, V(G)\setminus U)\cap \{6, 7, \ldots, m\}|\leq 1$ for each $i\in\{1, 3, 4, 6\}$. Since there are $m\geq {\rm max}\left\{\frac{k}{2}+8, 77\right\}$ colors used in $G$, we may assume that $c(u_2u)=m$ for some $u\in V(G)\setminus U$. In order to avoid a rainbow $K_2\cup 2P_3$, we have that $G[V(G)\setminus (U\cup\{u\})]$ is a monochromatic subgraph using color 4, contradicting to the assumption that $G$ contains no $ASMS(k)$.

Secondly, we consider the case when $H=2K_2\cup K_3$. By Corollary \ref{co:K2k3K2P4+K2P5P3P4}, there is a rainbow $K_2\cup K_3$ in $G$, say with vertex set $V=\{v_1, v_2, v_3, v_4, v_5\}$ such that $c(v_1v_2)=1$, $c(v_2v_3)=2$, $c(v_3v_1)=3$ and $c(v_4v_5)=4$. Then $C(V(G)\setminus V)\subseteq \{1, 2, 3, 4\}$ to avoid a rainbow $2K_2\cup K_3$. We first suppose that $C(\{v_4, v_5\}, V(G)\setminus V)\cap \{5, 6, \ldots, m\}\neq \emptyset$, say $c(v_4v)=5$ for some $v\in V(G)\setminus V$. Then $C(V(G)\setminus (V\cup\{v\}))\subseteq \{1, 2, 3\}$, and by Theorem \ref{th:2coloring} we have $C(V(G)\setminus (V\cup\{v\}))=\{1, 2, 3\}$. By Lemma \ref{le:conK3}, $G[V(G)\setminus (V\cup\{v\})]$ contains a rainbow $K_3$, say with vertex set $V'$. In order to avoid a rainbow $2K_2\cup K_3$, we have $C(V'\cup\{v\}\cup (V\setminus \{v_4\}), V(G)\setminus (V\cup V'\cup\{v\}))\subseteq \{1, 2, 3, 4, 5\}$. Since there are $m\geq 77> \binom{|V\cup V'\cup\{v\}|}{2}$ colors used in $G$, we may further assume that $c(v_4v')=m$ for some $v'\in V(G)\setminus (V\cup V'\cup\{v\})$. But then $c(vv')\notin \{1, 2, \ldots, m\}$, a contradiction. Therefore, $C(\{v_4,v_5\}, V(G)\setminus V)\cap \{5, 6, \ldots, m\}=\emptyset$. Since there are $m\geq {\rm max}\left\{\frac{k}{2}+8, 77\right\}$ colors used in $G$, we may assume that $|C(v_1, V(G)\setminus V)\cap \{5, 6, \ldots, m\}|\geq 2$, say $c(v_1v'')=5$ and $c(v_1v''')=6$ for some $v'', v'' \in V(G)\setminus V$. Then $c(v''v''')\in \{2, 4\}$. And we further have $C(V(G)\setminus (V\cup\{v'', v'''\}))\subseteq \{2, 4\}$. By Theorem \ref{th:2coloring}, there is an $ASMS(k)$, a contradiction.

Next, we consider the case when $H=2K_2\cup P^+_4$. By Corollary \ref{co:K2k3K2P4+K2P5P3P4}, there is a rainbow $K_2\cup P^+_4$ in $G$, say with vertex set $W=\{w_1, w_2, \ldots, w_7\}$ such that $c(w_1w_2)=1$, $c(w_2w_3)=2$, $c(w_3w_4)=3$, $c(w_3w_5)=4$ and $c(w_6w_7)=5$. Then $C(V(G)\setminus W)\subseteq \{1, 2, \ldots, 5\}$ for avoiding a rainbow $2K_2\cup P^+_4$. We first suppose that $C(\{w_6, w_7\}, V(G)\setminus W)\cap \{6, 7, \ldots, m\}\neq \emptyset$, say $c(w_6w)=6$ for some $w\in V(G)\setminus W$. Then $C(V(G)\setminus (W\cup\{w\}))\subseteq \{1, 2, 3, 4\}$, and by Theorem \ref{th:2coloring} we have $|C(V(G)\setminus (W\cup\{w\}))|\geq 3$. Let $c_1, c_2, c_3$ be three distinct colors used in $G[V(G)\setminus (W\cup\{w\})]$ and let $c(e_i)=c_i$ for $1\leq i\leq 3$, where $e_1$, $e_2$ and $e_3$ are three edges in $G[V(G)\setminus (W\cup\{w\})]$. Let $W'$ be the set of vertices incident with at least one of $e_1$, $e_2$ and $e_3$, so $|W'|\leq 6$. Let $A=V(G)\setminus(W\cup W'\cup \{w\})$. In order to avoid a rainbow $2K_2 \cup P^+_4$, we have $|C(u_i, A)\cap \{7, 8, \ldots, m\}|\leq 1$ for each $u_i\in W\setminus \{w_6\}$. Since there are $m\geq 77> \binom{|W\cup\{w\}|}{2} + |W||W'| + |W\setminus \{w_6\}|$ colors used in $G$, we may further assume that $c(w_6w')=m$ for some vertex $w'\in A$. Let $w''\in A\setminus \{w'\}$ and we will consider the color used on the edge $w'w''$. Note that $c(w'w'')\in\{1, 2, \ldots, 5\}$. If $c(w'w'')\in \{1, 2, 3, 4\}$, then $\{w_6, w_7, w, w', w''\}$ forms a rainbow $P^+_4$, which together with one edge within $\{w_1, w_2, \ldots, w_5\}$ and one edge within $W'$ forms a rainbow $2K_2\cup P^+_4$. If $c(w'w'')=5$, then $w'w''$, $w_6w$ and the rainbow $P^+_4$ within $\{w_1, w_2, \ldots, w_5\}$ form a rainbow $2K_2\cup P^+_4$, a contradiction.

Therefore, $C(\{w_6, w_7\}, V(G)\setminus W)\cap \{6, 7, \ldots, m\}=\emptyset$. Thus every edge with color in $\{6, 7, \ldots, m\}$ is incident with a vertex in $\{w_1, w_2, \ldots, w_5\}$. If $C(V(G)\setminus W)\cap \{3, 4\}\neq \emptyset$, say $c(ab)=3$ for some $a, b \in V(G)\setminus W$, then $|C(w_i, V(G)\setminus (W\cup \{a, b\}))\cap \{6, 7, \ldots, m\}|\leq 1$ for every $1\leq i\leq 5$, so the total number of colors used in $G$ is at most $5+\binom{|W|}{2}\leq 5+21<m$, a contradiction. Thus $C(V(G)\setminus W)=\{1, 2, 5\}$. Let $c(f_1)=1$, $c(f_2)=2$ and $c(f_3)=5$, where $f_1$, $f_2$ and $f_3$ are three edges in $G[V(G)\setminus W]$. Let $W''$ be the set of vertices incident with at least one of $f_1$, $f_2$ and $f_3$. In order to avoid a rainbow $2K_2\cup P^+_4$, we have $|C(w_i, V(G)\setminus (W\cup W''))\cap \{6, 7, \ldots, m\}|\leq 1$ for every $i\in \{2, 4, 5\}$.
For each $j\in \{1, 3\}$, if $|C(w_j, V(G)\setminus (W\cup W''))\cap \{6, 7, \ldots, m\}|\geq 3$, say $c(w_jx_1)=6$, $c(w_jx_2)=7$ and $c(w_jx_3)=8$, then for any vertex $y\in V(G) \setminus (W\cup W''\cup \{x_1, x_2, x_3\})$, we have $c(x_1y)\notin \{1, 2, \ldots, m\}$, a contradiction. Thus the total number of colors used in $G$ is at most $3+2\cdot 2+\binom{|W|}{2}\leq 7+21<m$, a contradiction.

Finally, we consider the case when $H=3K_2\cup K_{1,3}$. From the above argument, there is a rainbow $2K_2\cup P^+_4$ in $G$, say vertex set $B$ and color set $\{1, 2, \ldots, 6\}$. In order to avoid a rainbow $3K_2\cup K_{1,3}$, we have $C(V(G)\setminus B)\subseteq \{1, 2, \ldots, 6\}$ and $|C(b, V(G)\setminus B)\cap \{7, 8, \ldots, m\}|\leq 1$ for any vertex $b\in B$. Thus the total number of colors used in $G$ is at most $|B|+\binom{|B|}{2} \leq 9+36<m$, a contradiction.
\hfill$\blacksquare$
\vspace{0.1cm}

By Corollary \ref{co:K2k3K2P4+K2P5P3P4}, Fact \ref{fa:2}, Lemmas \ref{le:H} and \ref{le:H+}, Theorem \ref{th:discon} is true.

\section{Proof of Theorem \ref{th:main2}}
\label{sec:ch3}

We first state some known results which will be used in the proofs of our main results.

\noindent\begin{lemma}\label{le:rcol} {\normalfont (\cite{Gya})}
Every $m$-coloring of $K_n$ has a monochromatic connected subgraph with at least $\frac{n}{m-1}$ vertices.
\end{lemma}

\noindent\begin{lemma}\label{le:rcolKmn} {\normalfont (\cite{LiMP2})}
The order of the largest monochromatic component of an $m$-coloring of $K_{s,t}$ is at least $\frac{s+t}{m}$.
\end{lemma}

\noindent\begin{lemma}\label{le:rcolkcon} {\normalfont (\cite{LiMP2})}
Let $n, m, k$ be positive integers and $\epsilon >0$ satisfy $m\geq 3$ and $n \geq \frac{11(2+\epsilon)}{\epsilon}k^2m^2$. Then in any $m$-coloring of $K_n$, there is a $k$-connected monochromatic subgraph of order at least $\frac{n}{m-1} - \left(1+\frac{1}{m(m-2)}+\epsilon\right)k^2m$.
\end{lemma}

\noindent\begin{lemma}\label{le:mader} {\normalfont (\cite{Mad})}
Let $\alpha$ be a real number, and let $G$ be a graph with average degree $\alpha$. Then $G$ has an $\frac{\alpha}{4}$-connected subgraph $H$, and therefore $|V(H)|\geq \frac{\alpha}{4}$.
\end{lemma}

In the following, we will first consider the case when $k=1$ in Theorem \ref{th:main2}.

\noindent\begin{theorem}\label{th:k1}
For any graph $H\in \mathcal{G}\cup \mathcal{H}$, let $m\geq {\rm max}\left\{\frac{k}{2}+8, 77\right\}$ and $n\geq (m+11)(k-1)+7k {\rm ln} k+2m+3$. Then any rainbow $(P_3 \cup H)$-free coloring of $K_n$ using all the $m$ colors contains a connected monochromatic subgraph of order at least $\frac{n-|V(H)|}{|E(H)|}$.
\end{theorem}

\noindent{\bf Proof.} For any $H\in \mathcal{G}\cup \mathcal{H}$, let $G$ be a rainbow $(P_3 \cup H)$-free coloring of $K_n$ using all the $m$ colors, and suppose that $G$ contains no connected monochromatic subgraph of order at least $\frac{n-|V(H)|}{|E(H)|}$.

Since $H \subseteq \mathcal{G}\cup \mathcal{H}$, every rainbow $H$-free coloring of $K_n$ contains an almost spanning $k$-connected (and thus connected) monochromatic subgraph with order greater than $\frac{n-|V(H)|}{|E(H)|}$. Thus there is a rainbow $H$ in $G$, say with vertex set $V(H)$ and color set $C(H)$. Let $C'=C(G)\setminus C(H)$, where $C(G)=\{1, 2, \ldots, m\}$ is the set of colors used in $G$. If $|C(G-V(H)) \cap C'|\leq 1$ (that is, $|C(G-V(H))|\leq |E(H)|+1$), then there is a connected monochromatic subgraph of order at least $\frac{n-|V(H)|}{|E(H)|}$ by Lemma \ref{le:rcol}, a contradiction. Thus $|C(G-V(H)) \cap C'|\geq 2$. In order to avoid a rainbow $P_3\cup H$, for any two edges $e_1, e_2\in E(G-V(H))$ with $c(e_1)\neq c(e_2)$ and $c(e_1), c(e_2)\in C'$, we have that $e_1$ and $e_2$ are non-adjacent. Let $U=\{v\in V(G-V(H))\ |\ v$ is incident with an edge in $G-V(H)$ using color $m\}$ and $V= V(G-V(H))\setminus U$. Note that we have $U\neq \emptyset$ and $V\neq \emptyset$. Then $C(U, V)\subseteq C(H)$, and there is a connected monochromatic subgraph of order at least $\frac{|U|+|V|}{|E(H)|} = \frac{n-|V(H)|}{|E(H)|}$ by Lemma \ref{le:rcolKmn}, a contradiction.
\hfill$\blacksquare$
\vspace{0.1cm}

Next we will consider the case when $k\geq 2$ in Theorem \ref{th:main2}. Note that we have not tried to optimize the bound on $n$ and the order of the $k$-connected monochromatic subgraph in the following result, since our aim is to show that the largest $k$-connected monochromatic subgraph has order at least $cn$, where $c=c(H)$ is a constant, not depending on $n$, $m$ or $k$.

\noindent\begin{theorem}\label{th:k2}
For any graph $H\in \mathcal{G}\cup \mathcal{H}$, let $k\geq 2$, $m\geq {\rm max}\left\{\frac{k}{2}+8, 77\right\}$ and $n\geq 99k^2m^2$. Then any rainbow $(P_3 \cup H)$-free coloring of $K_n$ using all the $m$ colors contains a $k$-connected monochromatic subgraph of order at least $\frac{n-|V(H)|}{9|E(H)|}$.
\end{theorem}

\noindent{\bf Proof.} For any $H\in \mathcal{G}\cup \mathcal{H}$, let $G$ be a rainbow $(P_3 \cup H)$-free coloring of $K_n$ using all the $m$ colors, and suppose that $G$ contains no $k$-connected monochromatic subgraph of order at least $\frac{n-|V(H)|}{9|E(H)|}$. Since $H\in \mathcal{G}\cup \mathcal{H}$, we may further assume that there is a rainbow $H$ in $G$, otherwise there would be a $k$-connected almost spanning monochromatic subgraph of order greater than $\frac{n-|V(H)|}{9|E(H)|}$. Consider a rainbow copy of $H$ in $G$ with vertex set $V(H)$ and color set $C(H)=\{1, 2, \ldots, |E(H)|\}$. Let $C'=C(G)\setminus C(H)$, where $C(G)=\{1, 2, \ldots, m\}$ is the set of colors used in $G$.

We first have the following simple fact attained by setting $\epsilon = \frac{1}{4}$ in Lemma \ref{le:rcolkcon}.

\begin{fact}\label{fa:1} If $m\geq 3$, $k\geq 2$ and $n\geq 99k^2m^2$, then in any $m$-coloring of $K_n$, there is a $k$-connected monochromatic subgraph of order at least $\frac{n}{m-1} - \frac{19}{12}k^2m$.
\end{fact}

If $|C(G-V(H)) \cap C'|\leq 1$ (that is, $|C(G-V(H))|\leq |E(H)|+1$), then there is a $k$-connected monochromatic subgraph of order at least $\frac{n-|V(H)|}{|E(H)|} - \frac{19(|E(H)|+1)}{12} \cdot k^2\geq \frac{n-|V(H)|}{9|E(H)|}$ by Fact \ref{fa:1}, a contradiction. Thus $|C(G-V(H)) \cap C'|\geq 2$. In order to avoid a rainbow $P_3\cup H$, for any two edges $e_1, e_2\in E(G-V(H))$ with $c(e_1)\neq c(e_2)$ and $c(e_1), c(e_2)\in C'$, we have that $e_1$ and $e_2$ are non-adjacent. Let $U_i=\{v \in V(G)\setminus V(H))\ |\ v \mbox{ is incident with}$ an edge in $G-V(H)$ using color $i\}$ for every $i\in C'$, and let $U_{m+1}= V(G-V(H)) \setminus \left(\bigcup _{i\in C'} U_i\right)$. Recall that $|C(G-V(H)) \cap C'|\geq 2$, so at least two of $U_{|E(H)|+1}, \ldots, U_{m}$, $U_{m+1}$ are non-empty. Let $U_{i_1}, U_{i_2}, \ldots, U_{i_t}$ be all the non-empty sets of $U_{|E(H)|+1}, \ldots, U_{m}$, $U_{m+1}$, where $2 \leq t \leq m-|E(H)|+1$. Note that all the edges between these $t$ parts are colored by colors in $C(H)$, and $C(U_{j})\subseteq C(H)\cup\{j\}$ for all $j \in \{i_1, i_2, \ldots, i_t\}$. Then each $U_j$ ($j \in \{i_1, i_2, \ldots, i_t\}$) satisfies $|U_j|\leq \frac{n}{12}-1$, since otherwise there would be a $k$-connected monochromatic subgraph of order at least $\frac{|U_j|}{|E(H)|}-\frac{19(|E(H)|+1)}{12} \cdot k^2 \geq \frac{n}{12|E(H)|}-\frac{19(|E(H)|+1)}{12} \cdot k^2 \geq \frac{n-|V(H)|}{9|E(H)|}$ using Fact \ref{fa:1} with $n\geq 99k^2m^2$. If $t \leq 12$, then $|V(G-V(H))|\leq t\cdot (\frac{n}{12}-1)\leq 12\cdot (\frac{n}{12}-1) \leq n-12 < n-|V(H)|$, a contradiction. Thus $t\geq 13$.

We now choose a subset $L \subseteq \{1, 2, \ldots, t\}$ such that $A=\bigcup_{l \in L}U_{i_{l}}$ and $B=\bigcup_{h \in \{1, 2, \ldots, t\}\setminus L}$ $U_{i_{h}}$ satisfy

(1) $|A|\geq |B|$;

(2) $|A|-|B|$ is minimum subject to (1).

\begin{claim}\label{cl:2} $|A|-|B| \leq \frac{n-|V(H)|}{6}$.
\end{claim}

{\bf Proof.} If $|A|-|B| \geq \frac{n-|V(H)|}{6}+1$, then for any element $s \in L$, let $A'= A\setminus U_{i_s}$ and $B'= B\cup U_{i_s}$. Then we have $|A'|-|B'|= |A|-|U_{i_s}|-(|B|+|U_{i_s}|)= |A|-|B|-2|U_{i_s}|$. Thus $|A'|-|B'|\geq |A|-|B|-2\left(\frac{n}{12}-1\right)=|A|-|B|-\frac{n}{6}+2 \geq \frac{n-|V(H)|}{6}+1 -\frac{n}{6}+2 \geq 0$ and $|A'|-|B'| \leq |A|-|B|-2 <|A|-|B|$, contradicting to the choice of $L$.
\hfill$\square$
\vspace{0.1cm}

Since $0 \leq |A|-|B| \leq \frac{n-|V(H)|}{6}$ and $|A|+|B|=n-|V(H)|$, we have $|A||B|\geq \frac{5(n-|V(H)|)}{12} \cdot \frac{7(n-|V(H)|)}{12} = \frac{35(n-|V(H)|)^2}{144}$, that is, there are at least $\frac{35(n-|V(H)|)^2}{144}$ edges between $A$ and $B$. Since all the edges between $A$ and $B$ are colored by colors in $C(H)$, there are at least $\frac{35(n-|V(H)|)^2}{144|E(H)|}$ edges using a single color, say color 1. Consider the bipartite graph $H'$ with bipartition $(A, B)$ and $E(H')=\{uv\ |\ u\in A, v\in B, c(uv)=1\}$. The average degree of $H'$ is at least $\frac{2}{n-|V(H)|} \cdot \frac{35(n-|V(H)|)^2}{144|E(H)|} = \frac{35(n-|V(H)|)}{72|E(H)|}$. By Lemma \ref{le:mader}, there is a subgraph $H''\subseteq H'$ with connectivity $\frac{1}{4} \cdot \frac{35(n-|V(H)|)}{72|E(H)|} > \frac{n-|V(H)|}{9|E(H)|} > k$. Then $H''$ is a $k$-connected monochromatic subgraph of order at least $\frac{n-|V(H)|}{9|E(H)|}$, a contradiction.
\hfill$\blacksquare$
\vspace{0.1cm}

Since $n\gg |V(H)|$, if we let $c(H)=\frac{n}{10|E(H)|}$, then Theorem \ref{th:main2} holds by Theorems \ref{th:k1} and \ref{th:k2}.

\section{A classification of the set $\mathcal{B}$}
\label{sec:ch4*}

\noindent\begin{proposition}\label{prop:B1}
The set $\mathcal{B}$ is a subset of $2P_3$, $2K_2\cup K_{1,3}$ and their subgraphs of order at least 3.
\end{proposition}

\noindent {\bf Proof.} Let $K_{s,t}$ be a complete bipartite graph with $s$ vertices in one partite set $U$ and $t$ vertices in the other partite set $V$. We will consider three colorings of $K_{s,t}$ (see Fig. 4.1) constructed as follows. The coloring $F_1$ is an $m$-coloring of $K_{s,t}$ with $U=U_1\cup U_2\cup \cdots \cup U_m$ such that each of $U_1, U_2, \ldots, U_m$ contains about $\frac{s}{m}$ vertices. We color the edges such that $c(U_i, V)=i$ for every $1\leq i\leq m$. The coloring $F_2$ is an $m$-coloring of $K_{s,t}$ with $U=U_1\cup U_2\cup \{u\}$ such that each of $U_1, U_2$ contains about $\frac{s-1}{2}$ vertices. We color the edges such that $c(U_1, V)=1$, $c(U_2, V)=2$ and $C(u, V)=\{1, 2, \ldots, m\}$. The coloring $F_3$ is an $m$-coloring of $K_{s,t}$ with $U=U_3\cup U_4\cup \cdots \cup U_m$ and $V=V_3\cup V_4\cup \cdots \cup V_m$, such that each of $U_3, U_4, \ldots, U_m$ contains about $\frac{s}{m-2}$ vertices, and each of $V_3, V_4, \ldots, V_m$ contains about $\frac{t}{m-2}$ vertices. We color the edges such that $c(U_i, V_i)=i$ for every $3\leq i\leq m$, $c(\bigcup^{\alpha}_{j=3}U_{j}, \bigcup^{m}_{l=\alpha+1}V_{l})=c(\bigcup^{\alpha}_{j=3}V_{j}, \bigcup^{m}_{l=\alpha+1}U_{l})=1$, and all the remaining edges are colored with color 2, where $\alpha=\left\lfloor\frac{m-2}{2}\right\rfloor+2$. Note that the largest monochromatic $k$-connected subgraphs of $F_1$, $F_2$ and $F_3$ have order about $\frac{s}{m}+t$, $\frac{s-1}{2}+t$ and $\frac{s+t}{2}$, respectively. Thus every graph in $\mathcal{B}$ should be a rainbow subgraph of $F_1$, $F_2$ and $F_3$.

\begin{figure}[htb]
\begin{center}
\includegraphics {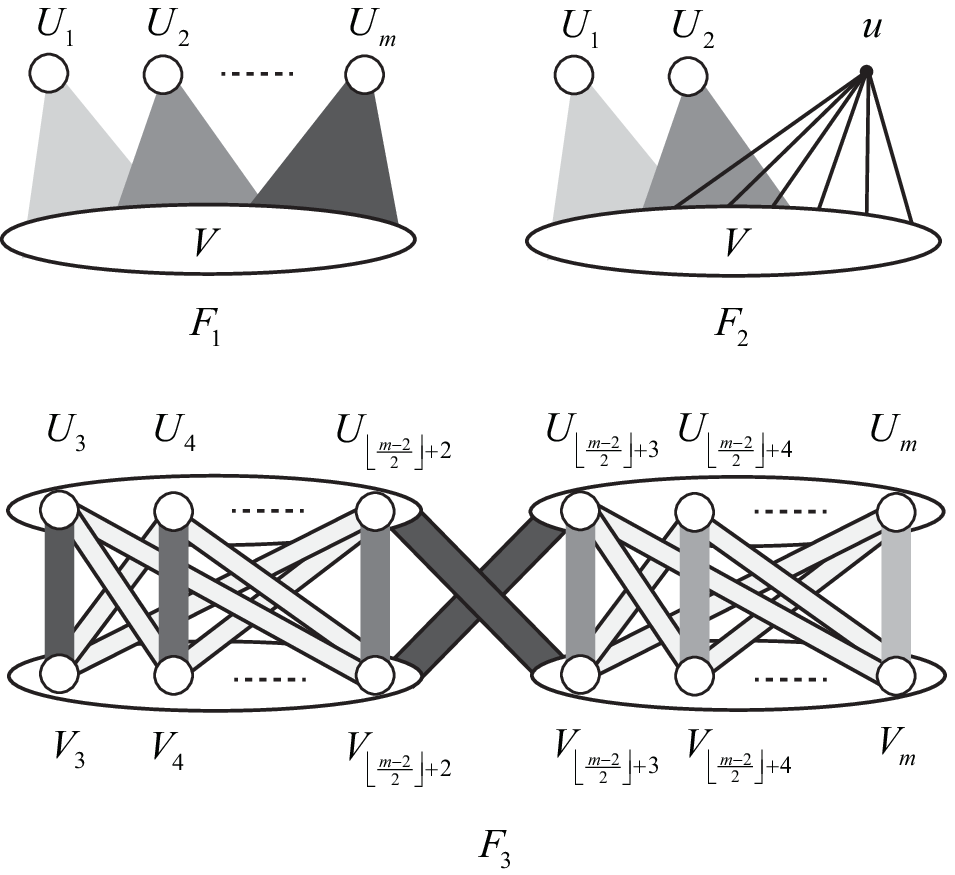}
\centerline{Fig. 4.1: The colorings $F_1$, $F_2$ and $F_3$.}
\end{center}
\end{figure}

For any graph $H\in \mathcal{B}$, since $F_1$ is rainbow $P_4$-free and $F_3$ is rainbow $K_{1,4}$-free, the longest path in $H$ has length at most 2, and the maximum degree of $H$ is at most 3. Thus each component of $H$ is one of $\{K_2, P_3, K_{1,3}\}$. Since $F_2$ is rainbow $4K_2$-free, $H$ has at most three components. Moreover, every rainbow $K_{1,3}$ in $F_3$ uses both color 1 and color 2, and every rainbow $P_3$ in $F_3$ uses at least one of color 1 and color 2. Thus $P_3\cup K_{1,3}\notin \mathcal{B}$. Furthermore, since $F_2$ is rainbow $K_2\cup 2P_3$-free, we have $K_2\cup 2P_3\notin \mathcal{B}$. From the above argument, we can deduce that $\mathcal{B}\subseteq\{P_3, K_{1,3}, 2K_2, K_2\cup P_3, K_2\cup K_{1,3}, 2P_3, 3K_2, 2K_2\cup P_3, 2K_2\cup K_{1,3}\}$. The result follows.
\hfill$\blacksquare$
\vspace{0.1cm}

Next, we character the structures of complete bipartite graphs without rainbow $K_{1,3}$ using a method similar to that used in our recent paper \cite{LiWL} with Liu.

\noindent\begin{theorem}\label{th:K13}
Let $K_{s,t}$, $s\geq t \geq 3$, be edge colored such that it contains no rainbow $K_{1,3}$.
Then, after renumbering the colors, one of the following holds:

{\rm (a)} at most four colors are used;

{\rm (b)} one partite set $U$ can be partitioned into $m-1$ non-empty parts $U_2, U_3, \ldots, U_m$, and the other partite set $V$ can be partitioned into $m-1$ non-empty parts $V_2, V_3, \ldots, V_m$, such that $C(U_i, V_i)\subseteq\{1, i\}$ for each $i\in \{2, 3, \ldots, m\}$, and every other edge is colored by color 1, where $m$ is the number of colors used in the coloring.
\end{theorem}

\noindent{\bf Proof.} Let $G$ be a rainbow $K_{1,3}$-free coloring of $K_{s,t}$ with partite sets $U=\{u_1, u_2, \ldots, u_s\}$ and $V=\{v_1, v_2, \ldots, v_t\}$. For a contradiction, suppose that there are $m\geq 5$ colors used in $G$ and (b) does not hold.

We call two adjacent edges of distinct colors a "$\wedge$". Note that every color appears on some $\wedge$. We claim that there are two $\wedge$s with four distinct colors (possibly such two $\wedge$s have common vertices, but no common edges). Indeed, since (b) does not hold, we may assume that colors 2 and 3 form a $\wedge$. In order to avoid two $\wedge$s with four distinct colors, each pair of colors in $C(G)\setminus \{2, 3\}$ cannot form a $\wedge$. Let $c(u'v')=c_1$, $c(u''v'')=c_2$, where $c_1, c_2\in C(G)\setminus \{2, 3\}$, $c_1 \neq c_2$, $u', u'' \in U$ and $v', v'' \in V$. For avoiding two required $\wedge$s, we have $c(u'v'')=c(u''v')=2$ or 3, say 2. For any color $c \in C(G)\setminus \{2, c_1, c_2\}$, there is no edge using color $c$ incident with $u'$, $u''$, $v'$ or $v''$ for avoiding two required $\wedge$s, so for every edge $uv$ with $c(uv)=c$, we may assume that $u \in U\setminus\{u', u''\}$ and $v \in V\setminus \{v', v''\}$. Then $c(u, \{v', v''\})=c(v, \{u', u''\})=2$. Thus, (b) holds if we exchange color 1 and color 2, a contradiction. Thus there are two $\wedge$s with four distinct colors. We may further consider four types of such two $\wedge$s (see Fig. 4.2).

\begin{figure}[htb]
\begin{center}
\includegraphics [width=0.85\textwidth, height=0.2\textwidth]{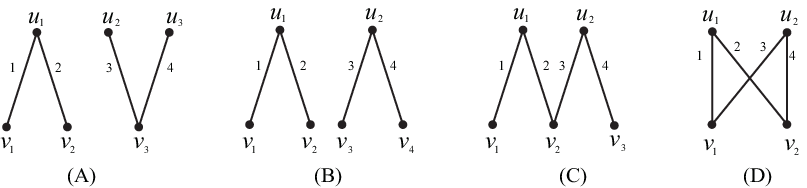}
\centerline{Fig. 4.2: Four types of $\wedge$s.}
\end{center}
\end{figure}

Firstly, if there are two $\wedge$s of type (A), then $c(u_1v_3)\notin\{1, 2, \ldots, m\}$, a contradiction. Thus there is no two $\wedge$s of type (A). Secondly, we consider type (B). For avoiding two $\wedge$s of type (A) and since $G$ is a rainbow $K_{1,3}$-free coloring using $m\geq 5$ colors, we may assume that $c(u_3v_5)=5$. But then $c(u_1v_5)\notin\{1, 2, \ldots, m\}$, a contradiction. Thus there is no two $\wedge$s of type (B), and there is no two $\wedge$s of type (C) for the same reason. Finally, we consider type (D). For avoiding two $\wedge$s of type (C) and since there are $m\geq 5$ colors used in $G$, we may assume that $c(u_3v_3)=5$. But then $c(u_1v_3)\notin\{1, 2, \ldots, m\}$ for avoiding a rainbow $K_{1,3}$ and two $\wedge$s of type (A), a contradiction. The result follows.
\hfill$\blacksquare$

\noindent\begin{corollary}\label{co:K13}
Given integers $k\geq 1$, $m\geq k+4$ and ${\rm min}\left\{s, t\right\}\geq m-1$, there is a spanning $k$-connected monochromatic subgraph in any rainbow $K_{1,3}$-free coloring of $K_{s,t}$ using all the $m$ colors.
\end{corollary}

\noindent{\bf Proof.} Let $G$ be a rainbow $K_{1,3}$-free coloring of $K_{s,t}$ with partite sets $U$ and $V$ using all the $m$ colors. Since $m\geq k+4\geq 5$ and by Theorem \ref{th:K13}, there is a partition $U_2, U_3, \ldots, U_m$ of $U$ and a partition $V_2, V_3, \ldots, V_m$ of $V$, such that $C(U_i, V_i)\subseteq\{1, i\}$ for each $i\in \{2, 3, \ldots, m\}$, and every other edge is colored by color 1.

Let $H$ be the subgraph induced by color 1. Then $|V(H)|=s+t$. For any vertex set $W$ of $k-1$ vertices $w_1, w_2, \ldots, w_{k-1} \in V(G)$, let $U'_i=U_i\setminus W$ and $V'_i=V_i\setminus W$ for every $2\leq i \leq m$. Since $m\geq k+4$, there are at least three distinct indexes $i_1, i_2, i_3\in \{2, 3, \ldots, m\}$ and three distinct indexes $j_1, j_2, j_3\in \{2, 3, \ldots, m\}$, such that $U'_{\alpha}=U_{\alpha}\neq\emptyset$ and $V'_{\beta}=V_{\beta}\neq\emptyset$ for each $\alpha \in \{i_1, i_2, i_3\}$ and each $\beta \in \{j_1, j_2, j_3\}$. Thus $H-W$ is connected, so $H$ is $k$-connected. The result follows.
\hfill$\blacksquare$

\noindent\begin{lemma}\label{le:2P3}
Given integers $k\geq 1$, $m\geq {\rm max}\left\{k+4, 7\right\}$ and ${\rm min}\left\{s, t\right\}\geq m$, there is a $k$-connected monochromatic subgraph of order at least $s+t-5$ in any rainbow $2P_3$-free coloring of $K_{s,t}$ using all the $m$ colors.
\end{lemma}

\noindent{\bf Proof.} For a contradiction, suppose that $G$ is a rainbow $2P_3$-free coloring of $K_{s,t}$ with partite sets $U=\{u_1, u_2, \ldots, u_s\}$ and $V=\{v_1, v_2, \ldots, v_t\}$ using all the $m$ colors, and $G$ contains no $k$-connected monochromatic subgraph of order at least $s+t-5$. By Corollary \ref{co:K13}, there is a rainbow $K_{1,3}$ in $G$, say $c(u_1v_i)=i$ for $1\leq i\leq 3$.

In order to avoid a rainbow $2P_3$, we have $|C(v_i, U\setminus \{u_1\})\cap \{4, 5, \ldots, m\}|\leq 1$ for every $1\leq i\leq 3$. If $|C(u_1, V\setminus \{v_1, v_2, v_3\})\cap \{4, 5, \ldots, m\}|\geq 1$, say $c(u_1v_4)=4$, then $|C(V(G)\setminus\{u_1, v_1, v_2, v_3, v_4\})|=1$ for avoiding a rainbow $2P_3$, which implies that there is a $k$-connected monochromatic subgraph of order at least $s+t-5$, a contradiction. Thus $C(u_1, V\setminus \{v_1, v_2, v_3\})\cap \{4, 5, \ldots, m\}=\emptyset$. Since there are $m\geq 7$ colors used in $G$, we may assume that $c(u_2v_4)=m$. Then $c(u_2, V\setminus \{v_1, v_2, v_3\}) = c(v_4, U\setminus \{u_1\})=m$ and thus $C(V(G)\setminus\{u_1, v_1, v_2, v_3\})=\{m\}$, resulting in a $k$-connected monochromatic subgraph of order at least $s+t-4$, a contradiction.
\hfill$\blacksquare$

\noindent\begin{lemma}\label{le:2K2K13}
Given integers $k\geq 1$, $m\geq {\rm max}\left\{k+4, 16\right\}$ and ${\rm min}\left\{s, t\right\}\geq m+1$, there is a $k$-connected monochromatic subgraph of order at least $s+t-6$ in any rainbow $2K_2\cup K_{1,3}$-free coloring of $K_{s,t}$ using all the $m$ colors.
\end{lemma}

\noindent{\bf Proof.}
Let $G$ be a rainbow $2K_2\cup K_{1,3}$-free coloring of $K_{s,t}$ with partite sets $U=\{u_1, u_2, \ldots,$ $u_s\}$ and $V=\{v_1, v_2, \ldots, v_t\}$ using all the $m$ colors. For a contradiction, suppose that $G$ contains no $k$-connected monochromatic subgraph of order at least $s+t-6$. By Corollary \ref{co:K13}, there is a rainbow $K_{1,3}$ in $G$, say $c(u_1v_i)=i$ for $1\leq i\leq 3$.

If $|C(u_1, V\setminus \{v_1, v_2, v_3\})\cap \{4, 5, \ldots, m\}|\geq 2$, say $c(u_1v_4)=4$ and $c(u_1v_5)=5$, then $|C(V(G)\setminus\{u_1, v_1, v_2, \ldots, v_5\})|=1$ for avoiding a rainbow $2K_2\cup K_{1,3}$, which implies that there is a $k$-connected monochromatic subgraph of order at least $s+t-6$, a contradiction. Thus $|C(u_1, V\setminus \{v_1, v_2, v_3\})\cap \{4, 5, \ldots, m\}|\leq 1$. For any $i\in\{1, 2, 3\}$, if $|C(v_i, U\setminus \{u_1\})\cap \{4, 5, \ldots, m\}|\geq 4$, say $c(v_iu_j)=j+2$ for every $j\in \{2, 3, 4, 5\}$, then $c(u_6v_4)\notin \{1, 2, \ldots, m\}$ for avoiding a rainbow $2K_2\cup K_{1,3}$, a contradiction. Thus $|C(v_i, U\setminus \{u_1\})\cap \{4, 5, \ldots, m\}|\leq 3$ for each $i\in\{1, 2, 3\}$.

Since there are $m\geq 16$ colors used in $G$, we have $|C(V(G)\setminus \{u_1, v_1, v_2, v_3\})\cap \{4, 5, \ldots, m\}|$ $\geq 3$. In order to avoid a rainbow $2K_2\cup K_{1,3}$, we may assume that there is a rainbow $K_{1,3}$ using colors $m-2$, $m-1$ and $m$ in $G-\{u_1, v_1, v_2, v_3\}$, say with vertex set $W$. Then for any edge $e\in E(G-W\cup \{u_1, v_1, v_2, v_3\})$, we have $c(e)\notin \{1, 2, \ldots, m\}$, a contradiction.
\hfill$\blacksquare$
\vspace{0.1cm}

By Proposition \ref{prop:B1}, Lemmas \ref{le:2P3} and \ref{le:2K2K13}, Theorem \ref{th:bipartite} is true.

\section{Large $k$-connected 2-colored subgraphs in Gallai-3-coloring}
\label{sec:ch4}

In this section, we will prove that Question \ref{ques-3} is true for $1\leq k \leq 3$, and we will give a counterexample to show that this question is false when $k=4t$, where $t$ is an positive integer. For this purpose, we will use the following structural result concerning Gallai-colorings given by Gallai \cite{Gallai}.

\noindent\begin{theorem}\label{th:Gallai} {\normalfont (\cite{Gallai})}
In any Gallai-coloring of a complete graph, there exists a partition $V_1, V_2, \ldots, V_l$ {\rm (}$l \geq 2${\rm )} of the vertices such that, there are in total at most two colors between the parts, and between every pair of parts there is only one color on the edges.
\end{theorem}

We now state and prove two lemmas that imply Question \ref{ques-3} is true for $1 \leq k \leq 3$ immediately.

\noindent\begin{lemma}\label{le:3/22}
For $n\geq 7$, every Gallai-3-coloring of $K_n$ contains a 2-connected subgraph of order $n$ using at most two colors.
\end{lemma}

\noindent{\bf Proof.}
Let $G$ be a Gallai-3-coloring of $K_n$ and suppose for a contradiction that $G$ contains no 2-connected spanning subgraph using at most two colors. We first prove the following claim concerning every 3-coloring of $K_n$.

\begin{claim}\label{cl:3/22}
For $n\geq 7$, there is a 2-connected subgraph of order $n-1$ using at most two colors in every 3-coloring of $K_n$.
\end{claim}

{\bf Proof.} For a contradiction, suppose that $F$ is a 3-coloring of $K_n$ using red, blue and green, which contains no 2-connected subgraph of order $n-1$ using at most two colors.

By Theorem \ref{th:2coloring}, there exists a 2-connected subgraph $R$ of order at least $n-2$ using either $\{$green$\}$ or $\{$red, blue$\}$. If $|V(R)|\geq n-1$, then $R$ is a desired subgraph so we may assume that the order of $R$ is exactly $n-2$. Let $V(F)\setminus V(R)=\{u,v\}$. First we consider the case that $R$ is a monochromatic subgraph colored by green. Since $|V(R)|=n-2\geq 5$, there are at least two edges between $u$ and $V(R)$ using a single color. Then $R$ and $u$ form a 2-connected subgraph of order $n-1$ using at most two colors, a contradiction. Therefore, we may further assume that $R$ is colored by red and blue. Then there is at most one edge using red or blue between $u$ (resp., $v$) and $R$. Hence, there is a subset $R'\subset V(R)$ with $|R'|=n-4 \geq 3$ such that $c(\{u,v\}, R')$ is green. Then there is a monochromatic complete bipartite graph $K_{2,n-4}$ colored by green, which is 2-connected. Let $V(R)\setminus R' = \{x, y\}$. Then we can find a 2-connected subgraph using at most two colors with vertex set $R' \cup \{u,v,x\}$ by a similar argument to the first case, a contradiction.
\hfill$\square$
\vspace{0.1cm}

By Claim \ref{cl:3/22}, $G$ contains a 2-connected subgraph $H$ of order $n-1$ using at most two colors, say red and blue. Let $V(G)\setminus V(H)= \{v\}$. Then there is at most one edge using red or blue between $v$ and $H$. First suppose that $c(vu)$ is red and $c(v, V(H)\setminus \{u\})$ is green for some $u \in V(H)$. Then to avoid a rainbow triangle, we have that $C(u, V(H)\setminus \{u\}) \subseteq \{\mbox{red, green}\}$. Then we obtain a $K_2 \vee \overline{K_{n-2}}$ using red and green, a 2-connected subgraph using at most two colors of order $n$. Thus we may further assume that $c(v, V(H))$ is green. Since any 2-coloring of complete graph contains a monochromatic spanning tree, there is a monochromatic spanning tree colored by either $\{$red$\}$ or $\{$blue, green$\}$ in $G[V(H)]$. In both cases, such a spanning tree together with vertex $v$ forms a 2-connected spanning subgraph of $G$ using at most two colors, a contradiction.
\hfill$\blacksquare$

\noindent\begin{lemma}\label{le:3/23}
For $n\geq 7$, every Gallai-3-coloring of $K_n$ contains a 3-connected subgraph of order at least $n-1$ using at most two colors.
\end{lemma}

\noindent{\bf Proof.}
Let $G$ be a Gallai-3-coloring of $K_n$, say using red, blue and green. Suppose that $G$ contains no 3-connected 2-colored subgraph of order at least $n-1$. By Lemma \ref{le:3/22}, $G$ contains a 2-connected spanning subgraph $H$ using at most two colors, say red and blue. Note that $H$ is not 3-connected, otherwise $H$ would be the desired subgraph. So we may further assume that $C=\{c_1, c_2\}$ is a cutset of $H$ and $(A, B)$ is a bipartition of the vertices of $H \setminus C$ such that $A$ (resp., $B$) is the union of vertices in components of $H \setminus C$. We may assume that $|A| \geq |B|$ without loss of generality. Moreover, all the edges are colored by green between $A$ and $B$, that is, $G[A\cup B]$ contains a monochromatic complete bipartite graph $H'$ of order $n-2$ using green.

If $|B|\geq 3$, then $H'$ is 3-connected. Furthermore, since $|A\cup B|=n-2\geq 5$, there are at least three edges between $c_1$ and $A\cup B$ using either $\{$red$\}$ or $\{$blue, green$\}$. If the former holds, then $G\setminus \{c_2\}$ contains a 3-connected subgraph using red and green of order $n-1$. And if the latter holds, then $G\setminus \{c_2\}$ contains a 3-connected subgraph using blue and green of order $n-1$. In both cases we can derive a contradiction. Thus $1 \leq |B| \leq 2$.

If $|B|=1$, say $B=\{b\}$, then $|A|=n-3\geq 4$. If $c(c_1b)$ (resp., $c(c_2b)$) is green, then $H\setminus \{c_2\}$ (resp., $H\setminus \{c_1\}$) is disconnected, contradicting the fact that $H$ is 2-connected. Thus $C(B, C)\subseteq \{\rm{red, blue}\}$. Without loss of generality, let $c_1b$ be colored by red. Then $C(c_1, A)\subseteq \{\rm{red, green}\}$ for avoiding a rainbow triangle. Likewise if $c_2b$ is colored by red, then $C(c_2, A)\subseteq \{\rm{red, green}\}$ holds once again. But then $G$ contains a 3-connected subgraph using red and green of order $n$, a contradiction. Thus $c_2b$ cannot be colored by red, so $c(c_2b)$ is blue. In order to avoid a rainbow triangle, we also have $C(c_2, A)\subseteq \{\rm{blue, green}\}$. Since $G[A]$ is a Gallai-3-coloring of $K_{|A|}$, there exists a color which spans a connected subgraph by Theorem \ref{th:Gallai}. If such a color is green or red, then $G\setminus \{c_2\}$ contains a 3-connected subgraph using green and red of order $n-1$, a contradiction. Thus this color is blue, but then $G\setminus \{c_1\}$ contains a 3-connected subgraph using blue and green of order $n-1$, a contradiction.

Hence $|B|=2$, say $B=\{b_1, b_2\}$. There exists at least one edge using red or blue between $B$ and $C$, say $c(c_1b_1)=$ red. Then $c(c_1, A)\subseteq \{\rm{red, green}\}$ to avoid a rainbow triangle. If $c(c_1b_2)$ is red or green, then $G\setminus \{c_2\}$ contains a 3-connected subgraph using red and green of order $n-1$, a contradiction. Thus $c(c_1b_2)$ is blue. In order to avoid a rainbow triangle, we have that $c(c_1, A)$ is green. But then then $H\setminus \{c_2\}$ is disconnected, contradicting the fact that $H$ is 2-connected. The proof of Lemma \ref{le:3/23} is complete.
\hfill$\blacksquare$
\vspace{0.1cm}

Finally, we shall construct a counterexample to show that Question \ref{ques-3} is false when $k= 4t$. We first state some results which will be used in our construction.

\noindent\begin{lemma}\label{le:EG} {\normalfont (\cite{ErGa*})}
Let $(d_1, d_2, \ldots, d_n)$ be a nonincreasing sequence of nonnegative integers. If $\sum^{n}_{i=1}d_i$ is even and
\begin{equation}
\sum_{i=1}^{k} d_{i} \leq k(k-1)+\sum_{i=k+1}^{n} \min \left\{k, d_{i}\right\}, \quad 1 \leq k \leq n
\end{equation}
then there is a simple graph with degree sequence $(d_1, d_2, \ldots, d_n)$.
\end{lemma}

\noindent\begin{corollary}\label{co:EG}
For all integers $t\geq 1$, let $d_1=d_2=\cdots=d_{2t}=2t$ and $d_{2t+1}=d_{2t+2}=\cdots=d_{4t}=2t-1$. Then there is a simple graph with degree sequence $(d_1, d_2, \ldots, d_{4t})$.
\end{corollary}

\noindent{\bf Proof.}
Since $\sum^{n}_{i=1}d_i=(2t)^2+2t(2t-1)$ is even, we need only to show that (1) holds by Lemma \ref{le:EG}. Let $a=2t$. If $1\leq k\leq a-1$, then
\begin{align*}
 k(k-1)+\sum_{i=k+1}^{2a} \min \left\{k, d_{i}\right\}- \sum_{i=1}^{k} d_{i}= & \ k(k-1)+ k(2a-k)-ak\\
   = & \ (a-1)k > 0.
\end{align*}
If $a\leq k\leq 2a$, then
\begin{align*}
 k(k-1)+\sum_{i=k+1}^{2a} \min \left\{k, d_{i}\right\}- \sum_{i=1}^{k} d_{i}= & \ k(k-1)+ (a-1)(2a-k)-a^2-(a-1)(k-a)\\
   = & \ k^2+2a^2-2ka-3a+k\\
   = & \ (k-a)^2+a^2-3a+k\\
   \geq & \ a^2-3a+a \geq 0.
\end{align*}
Thus (1) holds, and the result follows.
\hfill$\blacksquare$
\vspace{0.1cm}

By Corollary \ref{co:EG}, there is a graph $F$ of order $4t$ with degree sequence $(d_1, d_2, \ldots, d_{4t})$, where $d_1=d_2=\cdots=d_{2t}=2t$ and $d_{2t+1}=d_{2t+2}=\cdots=d_{4t}=2t-1$. We color all the edges of $F$ using color 1, and add an edge with color 2 between each pair of vertices that are not adjacent in $F$. Now we obtain a coloring $H$ of $K_{4t}$ using colors 1 and 2. Let $H'$ and $H''$ be the subgraphs of $H$ induced by color 1 and color 2, respectively. Without lose of generality, let $\{v_1, v_2, \ldots, v_{2t}\}$ and $\{v_{2t+1}, v_{2t+2}, \ldots, v_{4t}\}$ be the sets of vertices with degree $2t$ and $2t-1$ in $H'$, respectively. Then $\{v_1, v_2, \ldots, v_{2t}\}$ and $\{v_{2t+1}, v_{2t+2}, \ldots, v_{4t}\}$ are the sets of vertices with degree $2t-1$ and $2t$ in $H''$, respectively. Let $k=4t$, and then $\left\lfloor\frac{k-1}{2}\right\rfloor=2t-1$. Let $G$ be a coloring of $K_{n}$ with $V(G)=V_1\cup V_2\cup V_3$, where $|V_1|=n-6t$, $|V_2|=4t$ and $|V_3|=2t$. We color the edges of $G$ such that $G[V_2]=H$, $c(\{v_1, v_2, \ldots, v_{2t}\}, V_1)=1$, $c(\{v_{2t+1}, v_{2t+2}, \ldots, v_{4t}\}, V_1)=2$, and all the remaining edges are colored by color 3 (see Fig. 5.1 for an example).

\begin{figure}[htb]
\begin{center}
\includegraphics [width=0.4\textwidth, height=0.32\textwidth]{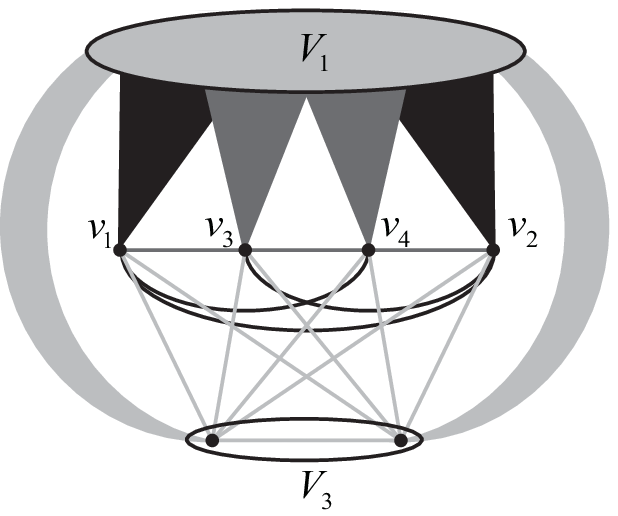}
\centerline{Fig. 5.1: A Gallai-3-coloring of $K_{n}$ without 4-connected 2-colored subgraph of order $n-1$.}
\end{center}
\end{figure}

We first show that $G$ is a Gallai-coloring. Suppose that there is a rainbow $K_3$ with vertex set $W$. Note that for any vertex $v\in V_3$, $c(v, V(G)\setminus \{v\})=3$. Thus $W\cap V_3=\emptyset$. Since $|C(V_1)|=1$ and $|C(V_2)|=2$, we have $|W\cap V_1|\leq 2$ and $|W\cap V_2|\leq 2$. Moreover, since $|C(V_2\cup \{u\})|=2$ for any vertex $u\in V_1$, we have $|W\cap V_1|= 2$ and $|W\cap V_2|= 1$. Since for any $v\in V_2$, $|C(v, V_1)|=1$, we also have $|W\cap V_2|= 2$ and $|W\cap V_1|= 1$, a contradiction. Thus $G$ is a Gallai-coloring.

Suppose that $G$ contains a $k$-connected subgraph $S$ of order at least $n- \left\lfloor \frac{k-1}{2} \right\rfloor=n-2t+1$ using at most two colors. Since $c(V_3, V_1\cup V_2)=3$ and $|V_3|=2t$, we have $3\in C(S)$. If $C(S)=\{1,3\}$, then $v_{2t+1}, v_{2t+2}, \ldots, v_{4t} \notin V(S)$ since there are exactly $4t-1$ edges using color 1 or color 3 between $x$ and $V(G)\setminus \{x\}$ for every $x\in \{v_{2t+1}, v_{2t+2}, \ldots, v_{4t}\}$. But then $|V(S)|\leq n-2t<n-2t+1$, a contradiction. Similarly, if $C(S)=\{2,3\}$, then we can also derive a contradiction. Thus $G$ contains no $k$-connected subgraph of order at least $n- \left\lfloor \frac{k-1}{2} \right\rfloor$ using at most two colors.

Although Question \ref{ques-3} is false when $k= 4t$, we still believe that the following conjecture is true.

\noindent\begin{conjecture}\label{conj-2}
For integers $k \geq 1$ and $n$ sufficient large, if $k\neq 4t$ for every integer $t\geq 1$, then every Gallai-3-coloring of $K_n$ contains a $k$-connected subgraph of order at least $n-\left\lfloor\frac{k-1}{2}\right\rfloor$ using at most two colors, and if $k= 4t$ for some integer $t\geq 1$, then every Gallai-3-coloring of $K_n$ contains a $k$-connected subgraph of order at least $n-\frac{k}{2}$ using at most two colors.
\end{conjecture}

\section{Concluding Remarks}
\label{sec:ch5}

In this paper, we mainly focus on the large $k$-connected monochromatic subgraph. Instead of finding a $k$-connected monochromatic subgraph, it will also be interesting to look for a long monochromatic path or cycle in colored complete graph. Using a result by Erd\H{o}s and Gallai \cite{ErGa} that for an integer $k \geq 2$ and a graph $G$ on $n$ vertices with $|E(G)|> \frac{k-1}{2}n$ there is a path $P_{k+1}$ in $G$, we prove the following result.

\noindent\begin{proposition}\label{prop:5}
Let $G$ be an $m$-coloring of $K_n$. For any non-negative integers $a_1, a_2, \ldots,$ $a_m$ with $\sum^{m}_{i=1}a_i \leq n+2m-2$, $G$ contains a monochromatic copy of $P_{a_i}$ in color $i$ for some $i \in [m]$.
\end{proposition}

\noindent{\bf Proof.} If $\mbox{min}_{1\leq i\leq m}\{a_1, a_2, \ldots, a_m\} \leq 2$, then $G$ contains a monochromatic $P_2$ clearly. Thus we may assume that $\mbox{min}_{1\leq i\leq m}\{a_1, a_2, \ldots, a_m\} \geq 3$. For $v\in V(G)$ and $i \in [m]$, let $d_{i}(v)$ denote the number of edges incident with $v$ using color $i$. Then for any $v \in V(G)$, we have $\sum^{m}_{i=1}d_i(v) =n-1 \geq \sum^{m}_{i=1}a_i - 2m + 1$. Let $\overline{d}_i(G)=\frac{1}{n}\sum_{v\in V(G)}d_i(v)$ for $1 \leq i \leq m$.
Then
\begin{align*}
 \sum^{m}_{i=1}\overline{d}_i(G) =\ & \sum^{m}_{i=1}\left(\frac{1}{n}\sum_{v\in V(G)}d_i(v)\right)=\ \frac{1}{n}\sum_{v\in V(G)}\left(\sum^{m}_{i=1}d_i(v)\right) \\
    \geq &\ \frac{1}{n}\sum_{v\in V(G)}\left(\sum^{m}_{i=1}a_i - 2m + 1\right) =\ \sum^{m}_{i=1}a_i - 2m + 1,
\end{align*}
and thus there exists an $i \in [m]$ such that $\overline{d}_i(G) > a_i - 2$. Then $G$ contains $\frac{1}{2}\sum_{v\in V(G)}d_i(v)= \frac{n}{2}\overline{d}_i(G) > \frac{a_i - 2}{2}n$
edges using color $i$. By Erd\H{o}s and Gallai's result mentioned above, $G$ contains a monochromatic copy of $P_{a_i}$ in color $i$.
\hfill$\blacksquare$

By setting $a_1=a_2=\cdots =a_m=\left\lfloor \frac{n+2m-2}{m} \right\rfloor$ in Proposition \ref{prop:5}, there is a monochromatic copy of $P_t$ with $t\geq \left\lfloor \frac{n+2m-2}{m} \right\rfloor =\left\lfloor \frac{n-2}{m} \right\rfloor + 2$ in any $m$-coloring of $K_n$. For monochromatic cycles, Kano and Li \cite{KaLi} showed that there is a monochromatic cycle of length at least $\left\lceil \frac{n}{m} \right\rceil$ in any $m$-coloring of $K_n$. It is natural to try to find the long paths or cycles in the coloring of $K_n$ in which we forbidden some rainbow subgraphs.

\end{document}